\documentclass[a4paper,10pt]{amsart}
\usepackage[matrix,arrow,tips,curve]{xy}
\usepackage[english]{babel}
\usepackage{amsmath}
\usepackage{amssymb}
\usepackage{pgf,tikz}
\usepackage{mathrsfs}
\usepackage{enumerate}
\usepackage{hyperref}
\usepackage{epstopdf}
\usepackage{multicol}
\usepackage{amssymb}
\usepackage{amsmath}
\usepackage{scalefnt}
\usepackage[T1]{fontenc}

\input{xy}
\xyoption{all}
\newtheorem{thm}{Theorem}[section]
\newtheorem{Lemma}[thm]{Lemma}

\newtheorem{Corollary}[thm]{Corollary}

\newtheorem*{thm*}{Theorem}

\theoremstyle{definition}

\newtheorem{Definition}[thm]{Definition}
\newtheorem{Remark}[thm]{Remark}

\newtheorem{Example}[thm]{Example}

\newcommand{\Hom}{\operatorname{Hom}}

\newcommand{\Pic}{\operatorname{Pic}}

\newcommand{\Proj}{\operatorname{Proj}}
\newcommand{\Sing}{\operatorname{Sing}}
\newcommand{\p}{\mathbb{P}}

\newcommand{\sF}{\mathscr{F}}
\newcommand{\sG}{\mathscr{G}}

\newcommand{\sI}{\mathscr{I}}

\newcommand{\sO}{\mathscr{O}}

\newcommand{\G}{\mathscr{G}}

\newcommand{\C}{\mathbb{C}}

\newcommand{\mm}{\mathfrak{m}}
\newcommand{\Z}{\mathbb{Z}}
\newcommand{\mO}{\mathcal{O}}

\DeclareMathOperator{\coker}{{coker}}
\newcommand{\inext}{{\mathcal E}{\it xt}}
\DeclareMathOperator{\supp}{Supp}

\newenvironment{dem}{\noindent {\bf Proof:}}
{\hfill \framebox[7pt]{} \mbox{} \medskip}

\begin{document}
\title{On holomorphic distributions on Fano Threefolds}
\author[Alana Cavalcante]{Alana Cavalcante}
\address{\sc Alana Cavalcante\\
UFOP\\
R. Trinta e Seis 115\\
35931-008, Jo\~ao Monlevade\\ Brazil}
\email{alana.decea@ufop.edu.br}

\author[Mauricio Corr\^ea]{Mauricio Corr\^ea}
\address{\sc Mauricio Corr\^ea\\
UFMG\\
Avenida Ant\^onio Carlos, 6627\\
30161-970, Belo Horizonte\\ Brazil}
\email{mauriciojr@ufmg.br}

\author[Simone Marchesi]{Simone Marchesi}
\address{\sc Simone Marchesi\\
UNICAMP\\
Rua Sergio Buarque
de Holanda, 651 \\
13083-859, Distr. Bar\~ao Geraldo,  Campinas\\ Brazil}
\email{marchesi@ime.unicamp.br}
\dedicatory{\it To  Israel Vainsencher on the occasion of his 70th birthday}

\date{\today}
\keywords{Holomorphic distributions, Fano manifolds, split vector
bundles, stable vector bundles }

\begin{abstract}
This paper is devoted to the study of holomorphic distributions of dimension and codimension
one on smooth weighted projective complete intersection Fano manifolds threedimensional, with Picard number equal to one.
We study the relations between algebro-geometric properties of the singular set of singular holomorphic distributions and their associated sheaves.
We characterize either distributions whose tangent
sheaf or  conormal sheaf are arithmetically Cohen Macaulay (aCM)  on smooth weighted projective complete intersection Fano manifolds threefold.
We also prove  that a  codimension one  locally free distribution with trivial canonical bundle  on any Fano threefold, with Picard number equal to one, has a tangent sheaf which either splits or it is stable.
\end{abstract}
\maketitle

\setcounter{tocdepth}{1}


\section{Introduction}

In complex manifolds, holomorphic distributions and foliations have been much studied (see \cite{ACM}, \cite{MOJ}, \cite{CJV}, \cite{GP}). An important problem is to analyse when the tangent and conormal sheaves split, as well as the properties of singular schemes of distributions. With this properties make it possible to classificate of codimension one distributions and study the stability of the tangent sheaf.

The relation between the tangent sheaf of a holomorphic distribution and its singular locus has been recently of great interest. Indeed, L. Giraldo and A. J. Pan-Collantes showed in \cite{GP} that the tangent
sheaf of a foliation of codimension  one  on $\mathbb{P}^3$  splits if and only
if its singular scheme $Z$ is aCM, and more recently,  M. Corr\^ea, M. Jardim and
R. Vidal Martins extended this result in \cite{CJV}, showing that the tangent sheaf of a
codimension one locally free distribution on $\mathbb{P}^n$ splits as a sum
of line bundles if and only if its singular scheme is aCM.


The section \ref{tangent sheaf} is
devoted precisely to extend this result for the smooth weighted projective complete intersection Fano threefolds with Picard number one.
Such varieties have been classified by Iskovskih \cite{I1,I2} and Mukai \cite{Mu}. The \textit{index} of $X$ is a basic invariant of these manifolds.
This is the largest integer $\iota_X$ such that the canonical line bundle $K_X$ is divisible by $\iota_X$ in $\Pic(X).$

In \cite{K-O}, Kobayashi and Ochiai showed that the index $\iota_X$ is at
most $\dim(X) + 1$ and $\iota_X = \dim(X)+1,$
if and only if $X \simeq \mathbb{P}^n.$ Moreover, $\iota_X = \dim(X),$ if and only if
$X \simeq Q^n \subset \mathbb{P}^{n+1},$ where $Q^n $ is a smooth quadric.
In particular, when $X$ is a Fano threefold, its can have index $\iota_X = 4,3,2,1.$
Separating the varieties by the index, we prove the following result.


\begin{thm}
Let $\sF$ be a distribution of codimension one on a smooth weighted projective complete intersection Fano threefold $X,$ such that the tangent sheaf $T_{\sF}$ is locally free. If $Z=\Sing(\sF)$  is the singular scheme  of $\sF$, then:
\begin{enumerate}
\item  If $\iota_X = 3$ and $T_{\sF}$ either splits as a sum of line bundles or is a spinor bundle, then  $Z$ is arithmetically Buchsbaum, with $h^1(Q^3, I_{Z} (r-2)) = 1$ being the only nonzero intermediate cohomology for $H^i(I_{Z}).$

\noindent Conversely, if $Z$ is arithmetically Buchsbaum
with $h^1(Q^3, I_{Z} (r-2)) = 1$ being the only nonzero intermediate cohomology for $H^i(I_{Z})$ and $h^2(T_{\sF}(-2)) = h^2 (T_{\sF} (-1-c_1(T_{\sF}))) = 0,$ then $T_{\sF}$ either split or is a spinor bundle.
\item If $\iota_X = 2$ and $T_{\sF}$ has no intermediate cohomology, then $H^1(X, I_{Z} (r+t)) = 0$ for $t < -6$ and $t > 8.$
Conversely, if $H^1(X, I_{Z} (r+t)) = 0$ for $t < -6$ and $t > 8,$ and $H^2(X, T_{\sF}(t)) = 0$ for $t \leq 8$ and
$H^1(X, T_{\sF}(s))  = 0$ for $s \neq -t- \iota_X -c_1(T_{\sF}),$ then $T_{\sF}$ has no intermediate cohomology.
\item If $\iota_X = 1$ and $T_{\sF}$ has no intermediate cohomology, then $H^1(X, I_{Z} (r+t)) = 0$ for $t < -4$ and $t > 4.$
Conversely, if $H^1(X, I_{Z} (r+t)) = 0$ for $t < -4$ and $t > 4,$ and $H^2(X, T_{\sF}(t)) = 0$ for $t \leq 4$ and
$ H^1(X, T_{\sF}(s))  = 0$ for $s \neq -t- \iota_X -c_1(T_{\sF}),$ then $T_{\sF}$ has no intermediate cohomology.

\end{enumerate}
\end{thm}

In section \ref{cotangent sheaf} we will focus on the conormal sheaf $N^{\ast}_{\sF}$ of a foliation of dimension one on Fano threefold.
M. Corr\^ea, M. Jardim and R. Vidal Martins showed
in \cite{CJV} that the conormal sheaf $N^{\ast}_{\sF}$ of a foliation of dimension one
on $\mathbb{P}^n$ splits if and only if its singular scheme $Z$ is arithmetically Buchsbaum
with $h^1(\mathcal{I}_Z(d-1)) = 1$ being the only nonzero intermediate cohomology.
We manage to extend the mentioned result for any Fano threefolds and obtain the following result.

\begin{thm}
Let $\sF$ be a distribution of dimension one on a smooth weighted projective complete intersection Fano threefold $X$, with index $\iota_X$ and Picard number one.
If $Z=\Sing(\sF)$  is the singular scheme  of $\sF$, then:
\begin{enumerate}
\item If $N^{\ast}_{\sF}$ is arithmetically Cohen Macaulay,  then $Z$ is arithmetically Buchsbaum, with $h^1(X, I_{Z} (r)) =1$ being the only nonzero intermediate cohomology for $H^i(I_{Z}).$
\item If $Z$ is arithmetically Buchsbaum with $h^1(X, I_{Z} (r)) = 1$ being the only nonzero intermediate cohomology for $H^i(I_{Z})$ and $h^2(N^{\ast}_{\sF}) = h^2 (N^{\ast}_{\sF} (-c_1 (N^{\ast}_{\sF}) - \iota_X)) = 0$ and $\iota_X \in \{1,2,3\},$ then $N^{\ast}_{\sF}$ is arithmetically Cohen Macaulay.
\end{enumerate}
\end{thm}

We addrees the properties of their singular schemes of codimension one holomorphic distributions on a threefold in section \ref{connect}.
D. Cerveau and A. Lins Neto  in \cite{CL}  have classified   codimension one foliations on $\mathbb{P}^3$ with   trivial canonical bundle. In \cite{MOJ}, O. Calvo-Andrade, M. Corr\^ea and M. Jardim  have studied  codimension distributions by  analyzing the properties of their singular schemes and tangent sheaves.  In \cite{LPT}, F. Loray, J. Pereira  and F. Touzet have studied codimension one foliation,  with numerically trivial canonical bundle, on Fano threefolds with Picard number one.

In \cite{MOJ}, O. Calvo-Andrade, M. Corr\^ea and M. Jardim showed a cohomological
criterion for the connectedness of the singular scheme of codimension one distributions
on $\p^3$ \cite[Theorem 3.7]{MOJ}. We extend this result for the others Fano threefolds with Picard number one. More precisely, we prove the following result.

\begin{thm}\label{con}
Let $\sF$ be a distribution of dimension one on a smooth weighted projective complete intersection Fano threefold $X$, with index $\iota_X$ and Picard number one.
Denote by $C$ the component of pure dimension one of $Z=\Sing(\sF)$.
 If $h^2(T_\sF(-r))=0$ and $C \subset X, \; C  \neq \emptyset,$ then $Z=C$ is connected and   $T_\sF$ is locally free.
Conversely, for $r \neq \iota_{X}$, if $Z = C$  and it is connected, then $T_\sF$ is locally free and $h^2(T_\sF(-r))=0$.
\end{thm}

In section \ref{stability} we will prove that if $\sF$ is a codimension one distribution on Fano threefold with tangent sheaf locally free and $c_1(T_{\sF}) = 0,$
then $T_{\sF}$ is either split or stable.

In \cite{MOJ}, O. Calvo-Andrade, M. Corr\^ea and M. Jardim showed that if $\sF$ is a codimension one distribution on $\mathbb{P}^3$ with $T_{\sF}$ is locally free and $c_1(T_{\sF})=0$, then $T_{\sF}$ is either split or stable. We prove that this in fact holds for any Fano threefold with  $X$ Picard number $\rho(X) = 1.$

\begin{thm}
Let $\sF$ be a codimension one distribution on a  Fano threefold with  $X$ Picard number $\rho(X) = 1.$ If  $T_{\sF}$ is locally free and $c_1(T_{\sF})=0$,
then $T_{\sF}$ is either split or stable.
\end{thm}

\hfill

\subsection*{Acknowledgments}  We are grateful to  Marcos Jardim for interesting conversations.  The first named author is grateful to  Department of Exact and Applied Science-UFOP;
she is also grateful to UFMG and Unicamp for their  hospitality.
The second named author was partially supported  by CNPq and  CAPES;   he is grateful to University of Oxford for its hospitality, in special he is thankful  to Nigel Hitchin for  interesting mathematical conversations.
The third author was partially supported by Funda\c c\~ao de Amparo a Pesquisa do Estado
de S\~ao Paulo (FAPESP) grant 2017/03487-9, by CNPq grant 303075/2017-1 and by the Coordenação de Aperfei\c coamento de
Pessoal de Nível Superior - Brasil (CAPES) - Finance Code 001.

\hfill

\subsection*{Notation and Conventions}
We always work over the field $\mathbb{C}$ of complex numbers. Given a complex variety $X$, we  denote by $TX$ the tangent bundle $(\Omega_{X}^1)^*$ and
to simplify the notation, given a distribution $\sF$ let us write $Z := \Sing(\sF).$

Assume that the Picard Group of $X$ is $\mathbb{Z}.$ We will denote
$E(t) = E \otimes_{\mathcal{O}_X} \mathcal{O}_X(t)$ for $t \in \mathbb{Z}$ when $E$ is a vector bundle on $X,$ and
$\mathcal{O}_X(1)$ denotes the its ample generator.

If $F$ is a sheaf on $X,$ we denote by $h^i(X, F)$ the dimension of the complex vector space
$H^i(X, F).$

\begin{Remark} \label{formula} \cite[Proposition 1.10]{H2}
For any holomorphic vector bundle $E$ of rank $r,$
$$\bigwedge^{k} E \simeq \bigwedge^{r-k}E^{\ast} \otimes \det E.$$
In particular, if $E$ is a rank $2$ reflexive sheaf, then
$ E^{\ast} = E \otimes (\det E)^{\ast}. $
\end{Remark}


\section{Preliminaries} \label{section:preliminaries}


In this section we recall some basic facts of the theory. Throughout this section, unless otherwise noted, $X$ denotes a smooth weighted projective complete intersection Fano threefold with Picard number one.


\stepcounter{thm}
\subsection{Fano manifolds} \label{section:Fano3folds}

A compact complex manifold $X$ is \textit{Fano} if
its anticanonical line bundle $\mathcal{O}_{X}(-K_{X}) \simeq \bigwedge^n TX$ is ample.

\begin{Definition} \label{wci}
We say that $X\subset \mathbb{P}(a_0,...,a_N)$ is a smooth $n$-dimensional \textit{weighted complete intersection} in a weighted projective space, when $X$ is the scheme-theoretic zero locus of $c = N-n$ weighted homogeneous polynomials $f_1,...,f_c$ of degrees $d_1,...,d_c$.
\end{Definition}

Smooth weighted projective complete intersection Fano threefold of Picard number one group have been classified by Iskovskikh \cite{I1,I2} and Mukai \cite{Mu}. They are:

\begin{itemize}
  \item [(a)] the projective space $\mathbb{P}^{3}$;
  \item [(b)] a quadric  hypersurface $Q^3 = X_2\subset\mathbb{P}^{4}$;
  \item [(c)] a cubic hypersurface $X_3\subset\mathbb{P}^{4}$;
  \item [(d)] an intersection $X_{2,2}$ of two quadric hypersurfaces in $\mathbb{P}^{5}$;
  \item [(e)] a hypersurface of degree $4$ in the weighted projective space $X_4\subset\mathbb{P}(1,1,1,1,2)$;
  \item [(f)] a hypersurface of degree $6$ in the weighted projective space $X_6\subset \mathbb{P}(1,1,1,2,3)$;
  \item [(g)] an intersection $X_{2,3}$ of a quadric and a cubic in $\mathbb{P}^5$;
  \item [(h)] an intersection $X_{2,2,2}$ of three quadrics in $\mathbb{P}^6$;
  \item [(i)] an intersection of a quadratic cone and a hypersurface of degree $4$ in $\mathbb{P}(1,1,1,1,1,2)$.
\end{itemize}

A basic invariant of a Fano manifold is its index.

\begin{Definition}
The \textit{index} of $X$ is the maximal integer $\iota_{X} > 0$ dividing $-K_{X}$ in $\Pic (X),$ i.e.
$-K_{X} = \iota_{X} \cdot H,$ with $H$ ample.
\end{Definition}

\begin{thm} \cite{K-O} \label{kobaiachi}
Let $X$ be Fano, $\dim(X) = n.$ Then the index $\iota_{X}$ is at
most $n + 1;$ moreover, if $\iota_{X} = n + 1,$ then $X \simeq \mathbb{P}^n,$ and if $\iota_{X}  = n,$ then $X$
is a quadric hypersurface $Q^n \subset \mathbb{P}^{n+1}.$
\end{thm}

By Theorem \ref{kobaiachi}, a Fano threefold $X$ can have $\iota_{X} \in \{1,2,3,4\}.$ Then, $\iota_{X} =4$ implies $X = \mathbb{P}^3,$
while $\iota_{X} =3,$ implies that $X$ is a smooth quadric hypersurface $Q^3$ in $\mathbb{P}^4.$ In case $\iota_{X} = 2,$ the variety $X$ is called a del Pezzo threefold, while $\iota_{X} = 1,$ the variety $X$ is called a prime Fano threefold.

By \cite[Theorem 3.3.4]{dolgachev},
\stepcounter{thm}
\begin{equation}\label{canwci}
K_X\cong\mathcal{O}_X\left(\sum_{j=1}^cd_j-\sum_{i=0}^{N}a_i\right)
\end{equation}

In particular, when $X$ is Fano, its index is
\stepcounter{thm}
\begin{equation} \label{indice}
\iota_X := \sum_{i=0}^{N}a_i-\sum_{j=1}^cd_j.
\end{equation}

Let $S_t$ be the $t$-th graded part of $S/(f_1,...,f_c)$. By \cite[Lemma 7.1]{CR},
\stepcounter{thm}
\begin{equation}\label{owci}
H^i(X,\mathcal{O}_X(t))\cong
\left\lbrace\begin{array}{lll}
S_t & \rm{if} & i=0; \\
0 & \rm{if} & 1\leq i\leq n-1;\\
S_{-t+\iota_X} & \rm{if} & i  = n.
\end{array} \right.
\end{equation}

\medskip

Let $Y$ be an $n$-dimensional Fano manifold  with $\rho(Y)=1$, and denote by $\sO_Y(1)$ the ample generator of
$\Pic(Y)$. Let $X \in \big|\sO_Y(d) \big|$ be a smooth divisor. We have the following exact sequences:
\stepcounter{thm}
\begin{equation}\label{restriction2}
0 \ \to \ \Omega_Y^{q}(t-d) \ \to \ \Omega_Y^{q}(t) \ \to \ \Omega_Y^{q}(t)|_X \ \to \ 0,
\end{equation}
and
\stepcounter{thm}
\begin{equation}\label{restriction1}
0 \ \to \ \Omega_X^{q-1}(t-d) \ \to \ \Omega_Y^{q}(t)|_X \ \to \ \Omega_X^{q}(t) \ \to \ 0.
\end{equation}

\medskip

\stepcounter{thm}
\subsection{Cohomology of cotangent sheaf} \label{cohomologia}

In this subsection we address important facts about the cohomology of $X$
with coefficients in an analytic coherent sheaf.





Let $p,q$ and $t$ be integers, with $p$ and $q$
non-negative. Then,
\begin{equation*} \label{bottP}
h^p\big(\p^n,\Omega_{\p^n}^q(t)\big) =
\begin{cases}
\binom{t+n-q}{t}\binom{t-1}{q} & \text{ for } p=0, 0\le q\le n \text{ and } t>q,\\
1 & \text{ for } t=0 \text{ and } 0\le p=q\le n,\\
\binom{-t+q}{-t}\binom{-t-1}{n-q} & \text{ for } p=n, 0\le q\le n \text{ and } t<q-n,\\
0 & \text{otherwise.}
\end{cases}
\end{equation*}

Throughout the work, we refer to this formula as classical Bott's formula.

In \cite{S}, Snow showed a vanishing theorem for the cohomology of $\Omega^q(t)$ for quadric hypersurfaces $X=Q^n$ in $\p^{n+1}$ and for a Grassmann manifold $X= Gr(s,m)$ of $s$-dimensional subspaces of $\mathbb{C}^m.$

\begin{thm} \cite{S} [Bott's formula for Quadric]\label{BottQ}
Let $X$ be a nonsingular quadric hypersurface of dimension $n.$
\begin{enumerate}
  \item If $-n+q \leq k \leq q$ and  $k \neq 0$ and $k \neq -n+2q,$ then $H^p(X, \Omega^q(k)) = 0$ for all p;
  \item $H^p(X, \Omega^q) \neq 0$ if and only if $p=q;$
  \item $H^p(X, \Omega^q(-n+2q)) \neq 0$ if and only if $p = n-q;$
  \item If $k > q,$ then $H^p(X, \Omega^q(k)) \neq 0$ if and only if $p=0;$
  \item If $k< -n+q,$ then $H^p(X, \Omega^q(k)) \neq 0$ if and only if $p=n.$
\end{enumerate}
\end{thm}


In \cite{dolgachev}, Dolgachev generalized the classical Bott's formula for the cohomology of twisted
sheaves of differentials to the case of weighted projective spaces ($\overline{\Omega}^{q}_{\p}(t)$) as follows:


Let $\p=\p(a_0,\ldots, a_N)=\Proj \big(S(a_0,\ldots, a_N)\big)$ be a weighted projective space.

Consider the sheaves of $\sO_{\p}$-modules $\overline{\Omega}^{q}_{\p}(t)$ defined in \cite[Section 2.1.5]{dolgachev} for $q, t\in\mathbb{Z}$, $q\ge 0$.
If $U\subset \p$ denotes the smooth locus of $\p$, and $\sO_{U}(t)$ is the line bundle obtained by restricting
$\sO_{\p}(t)$ to $U$, then ${\overline{\Omega}^{q}_{\p}(t)}_{|U}=\Omega^q_{U}\otimes \sO_{U}(t)$.
The cohomology groups  $H^p\big(\p, \overline{\Omega}^{q}_{\p}(t)\big)$ are described by:

\begin{thm} \cite[Section 2.3.2]{dolgachev} \label{Bott dolgachev} Let $h^p\big(\p, \overline{\Omega}^{q}_{\p}(t)\big) = \dim H^p \big(\p, \overline{\Omega}^{q}_{\p}(t)\big).$ Then:
\begin{itemize}
	\item[-] $h^0\big(\p, \overline{\Omega}^{q}_{\p}(t)\big)=\sum_{i=0}^q \Big((-1)^{i+q} \sum_{\#J=i}\dim_{\mathbb{C}}\big(S_{t-a_J}\big)\Big)$, where $J\subset \{0, \ldots, N\}$ and $a_J:=\sum_{i\in J}a_i$;
	\item[-] $h^0\big(\p, \overline{\Omega}^{q}_{\p}(t)\big) = 0$ if $t <\min \{\sum_{j\in J}a_{i_j}\: | \: \#J = q\}$;
	\item[-] $h^p\big(\p, \overline{\Omega}^{q}_{\p}(t)\big)=0$ if $p\not\in \{0, q,N\}$;
	\item[-] $h^p\big(\p, \overline{\Omega}^{p}_{\p}(t)\big)=0$ if $t\neq 0$ and $p\notin\{0,N\}$.
\end{itemize}
In particular, if $q\geq 1$, then
\stepcounter{thm}
\begin{equation}\label{lcwps}
h^0(\mathbb{P},\Omega^q_{\mathbb{P}}(t)) = 0 \ \text{ for any } \ t\leq q.
\end{equation}
\end{thm}

Notice that $\p(a_0,\ldots, a_N) = \mathbb{P}^N$ is a projective space we obtain the classical Bott's formulas.

Now assume that $\p$ has only isolated singularities, let $d>0$ be such that $\sO_{\p}(d)$ is a line bundle generated by global sections, and $X\in \big|\sO_{\p}(d)\big|$ a smooth hypersurface. We will use the cohomology groups $H^p\big(\p, \overline{\Omega}^{q}_{\p}(t)\big)$ to compute some cohomology groups
$H^p\big(X,\Omega_{X}^q(t)\big)$. Note that $X$ is contained in the smooth locus of $\p$, so we have an exact sequence as in \eqref{restriction1}:
\stepcounter{thm}
\begin{equation}\label{restriction1_P}
0 \to \ \Omega_X^{q-1}(t-d) \ \to \ \overline{\Omega}^{q}_{\p}(t)|_X \ \to \ \Omega_X^{q}(t) \to  0.
\end{equation}
Tensoring the sequence
$$
0 \to\ \sO_{\p}(-d)  \ \to \ \sO_{\p}  \ \to \ \sO_{X} \ \to 0.
$$
with the sheaf $\overline{\Omega}^{q}_{\p}(t)$, and noting that $\overline{\Omega}^{q}_{\p}(t)\otimes \sO_{\p}(-d) \cong \overline{\Omega}^{q}_{\p}(t-d)$, we get an exact sequence  as in \eqref{restriction2}:
\stepcounter{thm}
\begin{equation}\label{restriction2_P}
0 \to \overline{\Omega}^{q}_{\p}(t-d) \ \to \ \overline{\Omega}^{q}_{\p}(t) \ \to \ \overline{\Omega}^{q}_{\p}(t)|_X \to 0.
\end{equation}

\medskip

Now, we will conclude this section by recalling the cohomological computations on weighted complete intersections $X.$ The next theorem in terms of cohomology of $X,$ is due to Flenner:

\begin{thm} \label{flenner} \cite[Satz 8.11]{flenner81} Let $X$ be a weighted complete intersection. Then,
\begin{enumerate}
	\item[-] $h^q(X,\Omega_{X}^q) =1 $ for $0\le q\le n$, $q\neq \frac{n}{2}$.
	\item[-] $h^p\big(X,\Omega_{X}^q(t)\big) = 0$ in the following cases
		\begin{itemize}
			\item[-] $0<p<n$, $p+q\neq n$ and either $p\neq q$ or $t\neq 0$;
			\item[-] $p+q > n$ and $t>q-p$;
			\item[-] $p+q < n$ and $t<q-p$.
		\end{itemize}
\end{enumerate}
\end{thm}

\medskip
\stepcounter{thm}
\subsection{aCM and aB schemes}

A closed subscheme $Y \subset X$ is \textit{arithmetically Cohen-Macaulay (aCM)} if its
homogeneous coordinate ring $S(Y) = k[x_0, \ldots, x_n]/I(Y)$ is a Cohen-Macaulay
ring.

Equivalently, $Y$ is aCM if $H_{\ast}^p(\mathcal{O}_Y) = 0$ for $1 \leq p \leq \dim(Y)-1$
and $H_{\ast}^1 (\mathcal{I}_Y) = 0$ (cf. \cite{CaH}).
From the long exact sequence of cohomology associated to the short exact sequence
$$
0 \rightarrow {\mathcal I}_Y \rightarrow \mathcal{O}_{X} \rightarrow \mathcal{O}_Y \rightarrow
0
$$
one also deduces that $Y$ is aCM if and only if $H_*^p({\mathcal
I}_Y)=0$ for $1\leq p\leq\dim\,(Y)$.

Similarly, a closed subscheme in $X$ is \textit{arithmetically
Buchsbaum (aB)} if its homogeneous coordinate ring is a Buchsbaum ring (see \cite{SV}).
Clearly, every aCM scheme is arithmetically Buchsbaum, but the
converse is not true: the union of two disjoint lines is
arithmetically Buchsbaum, but not aCM.

\medskip

\stepcounter{thm}
\subsection{Rank 2 vector bundles on quadrics }

Let $E$ be rank $2$ vector bundle on a smooth quadric $Q^n \subset \p^{n+1},$ with $n \geq 3.$
We set
$$R = \mathop\oplus_{t \geq 0} H^0(Q^n,\mathcal{O}_{Q^n}(t))\quad\text{and}\quad \mm = \mathop\oplus_{t>0} H^0(Q^n,\mathcal{O}_{Q^n}(t)),$$
and also
$$H^i_*(Q^n,E) = \mathop\oplus_{t \in \mathbb{Z}} H^i(Q^n,E(t))\quad \text{for }i=0,\ldots,n,$$
which are modules of finite length on the ring $R$. \\

\begin{Definition}
We say that $E$ is \emph{$k$-Buchsbaum}, with $k\geq0$,
 if for all integers $p$, $q$ such that $1\leq p\leq q-1$ and $3\leq q\leq n$ it holds
$$\mm^k \cdot H^p_*(Q',E\vert_{Q'}) = 0,$$
where $Q'$ is a general $q$-dimensional linear section of $Q^n$, i.e.\ $Q'$ is a quadric hypersurface cut out on $Q^n$ by a general linear space $L\subset \p^{n+1}$ of dimension $q+1$, that is $Q'=Q_n \cap L$.
\\
\end{Definition}

\begin{Remark}
Note that $E$ is $0$-Buchsbaum if and only if $E$ has no intermediate cohomology, i.e.\  $H^i(Q^n,E(t))=0$ for every $t\in \mathbb{Z}$ and $1\leq i\leq n-1$. Such a bundle is also called \emph{arithmetically Cohen-Macaulay}.
\\
Observe also that $E$ is $1$-Buchsbaum if and only if $E$ has every intermediate cohomology module with trivial structure over $R$. Such a bundle is also called  \emph{arithmetically Buchsbaum}.
\end{Remark}

\medskip

Now, we recall the definition and some properties of spinor bundles on
$Q^n.$ For more details see \cite{Ott1,Ott2}.

Let $S_k$ be spinor variety which parametrizes the family of
$(k-1)-$planes in $Q^{2k-1}$ or one of the two disjoint families
of $k-$planes in $Q^{2k}.$

We have $\dim(S_k) = (k(k+1))/2, \; \Pic(S_k) = \mathbb{Z}$ and $h^0 (S_k,\mathcal{O}(1)) = 2^k.$ Spinor varieties
are rational homogeneous manifolds of rank $1.$
When $n=2k-1$ is odd, consider for all $x \in Q^{2k-1}$ the variety
$$\{\mathbb{P}^{k-1} \in Gr (k-1,2k)/ x \in \mathbb{P}^{k-1} \subset Q^{2k-1}\}.$$
This variety is isomorphic to $S_{k-1}$ and we denote it by $(S_{k-1})_x.$ Then we have a natural embedding
$$(S_{k-1})_x \rightarrow S_k.$$
Considering the linear spaces spanned by these varieties, we have for all $x \in Q^{2k-1}$ a natural inclusion
$H^0 ((S_{k-1})_x, \mathcal{O}(1))^{\ast} \rightarrow H^0 (S_k, \mathcal{O}(1))^{\ast}$ and then an embedding
$$ s: Q^{2k-1} \rightarrow Gr(2^{k-1}-1,2^{k}-1),$$
in the Grassmannian of $(2^{k-1}-1)$-subspaces of $\mathbb{P}^{2^{k}-1}.$


\medskip
It is well known that $S_1 \simeq \mathbb{P}^1, \; S_2 \simeq \mathbb{P}^3.$
The embedding $s: Q^3 \rightarrow Gr(1,3)$ corresponds to a hyperplane section.

\begin{Definition}
Let $U$ be the universal bundle of the Grassmannian. We call $s^{\ast}U = S$
the \textit{spinor bundle} on $Q^{2k-1},$ whose rank is $2^{k-1}.$
\end{Definition}

\medskip
The spinor bundle $S$ on $Q^3$ is just the restriction of the universal sub-bundle on the $4$-dimensional quadric.
\medskip

\begin{Remark} \label{spinor glob gene}
The spinor bundle on $Q^3$ is a globally generated vector bundle of rank $2$, by construction.
\end{Remark}

The next theorem shows the existence of bundles that are indecomposable on quadric.

\begin{thm} \label{teo 2.1}\cite[Theorem 2.1]{Ott1}
The spinor bundle on $Q^3$ is stable.
\end{thm}

We have that $S$ is the unique indecomposable bundle of rank $2$ on $Q^3$ with
$c_1(S) = -1$ and $c_2(S) = 1,$ (see \cite{Bu}).

\medskip

\stepcounter{thm}
\subsection{Holomorphic Distributions}\label{holomorphic distributions}

Let us now recall facts about holomorphic distributions on complex projective varieties.
For more details about distributions and foliations see \cite{EK,CJV,CMS,MOJ}.

\begin{Definition} \label{distributionTF}
Let $X$ be a smooth complex manifold.
   A \textit{codimension $k$ distribution} $\mathscr{F} $ on $X$ is given by an exact sequence
\begin{equation}\label{eq:Dist}
\mathscr{F}:\  0  \longrightarrow T_{\sF} \stackrel{\phi}{ \longrightarrow} TX \stackrel{\pi}{ \longrightarrow} N_{\sF}  \longrightarrow 0,
\end{equation}
where $T_{\sF}$ is a coherent sheaf of rank $r_{\sF}:=\dim(X)-k$, 
and $N_{\sF}:= TX / \phi (T_{\sF})$ is a torsion free sheaf.

  The sheaves $T_{\sF}$ and $N_{\sF}$ are called the \textit{tangent} and the \textit{normal} sheaves of $\sF$, respectively.

  The \textit{singular set} of the distribution $\sF$ is defined by $\Sing(\sF) = \{x \in X | (N_{\sF})_x  \; \mbox {is not a free} \; \mathcal{O}_{X,x}-\mbox{module}\}.$

  The  \textit{conormal} sheaf  of $\sF$ is $N_{\sF}^*$.

\end{Definition}

A distribution $\sF$ is said to be {\it locally free} if $T_{\sF}$ is a locally free sheaf.

By definition, $\Sing(\sF)$ is the singular set of the sheaf $N_{\sF}.$ It is a closed
analytic subvariety of $X$ of codimension at least two, since by definition $N_{\sF}$ is torsion free.

\begin{Definition}
A \textit{foliation} is an integrable distribution, that is a distribution
$$ \mathscr{F}:\  0  \longrightarrow T_\sF \stackrel{\phi}{ \longrightarrow} TX \stackrel{\pi}{ \longrightarrow} N_{\sF}  \longrightarrow 0 $$
whose tangent sheaf is closed under the Lie bracket of vector fields, i.e. $[\phi(T_\sF),\phi(T_\sF)]\subset \phi(T_\sF)$.
\end{Definition}

Clearly, every distribution of codimension $\dim(X)-1$ is integrable.


\medskip


From now on, we will consider codimension one distributions on Fano threefold $X.$ Thus, sequence (\ref{eq:Dist}) simplifies to

\begin{equation}\label{eq:Dist Fano}
\mathscr{F}:\  0  \longrightarrow T_\sF \stackrel{\phi}{ \longrightarrow} TX \stackrel{\pi}{ \longrightarrow}
I_{Z/X}(r)  \longrightarrow 0,
\end{equation}
where $T_\sF$ is a rank $2$ reflexive sheaf and $r$ is an integer such that $r = c_1(TX) - c_1(T_{\sF})$. Observe that   $N_{\sF}=I_{Z/X}(r)$ where $Z$ is the singular scheme of $\sF$.

A codimension one distribution on $X$ can also be represented by a
section $\omega \in H^0 (\Omega^1_{X}(r))$, given by the dual of the morphism $\pi:TX \to N_{\sF}$. On the other hand, such section induces a sheaf map $\omega :\mathcal{O}_X \rightarrow \Omega^1_{X}(r)$. Taking duals, we get a cosection $$\omega^{\ast}:(\Omega^1_{X}(r))^{\ast}=TX(-r)\rightarrow \mathcal{O}_X$$
whose image is the ideal sheaf $I_{Z/X}$ of the singular scheme. The kernel of $\omega^{\ast}$ is the tangent sheaf $\sF$ of the distribution twisted by $\mathcal{O}(-r)$.

\begin{Remark}
From this point of view, the integrability condition is equivalent to $\omega \wedge d\omega=0.$
\end{Remark}


\section{Tangent sheaf vs. singular scheme} \label{tangent sheaf}


In this section, all the distributions that we will consider will be defined for Fano manifolds as described in Section \ref{section:Fano3folds}, of codimension one, and with singular set of codimension $\geq 2.$

For codimension one foliations in $\p^3,$ L. Giraldo and A. J. Pan Collantes showed in \cite{GP} that the tangent sheaf is locally free if and only if the singular scheme is a curve, i.e. a Cohen-Macaulay scheme
of pure dimension one, cf. [\cite{GP}, Theorem 3.2]. More recently, O. Calvo-Andrade, M. Corr\^ea, and M. Jardim obtained in \cite{MOJ} the following generalization of this theorem.

\begin{Lemma} \label{local-free-tg-sheaf}\cite[Lemma 2.1]{MOJ}
The tangent sheaf of a codimension one distribution is locally free if and only if its singular locus has pure codimension $2.$
\end{Lemma}

\medskip
We characterize when the tangent sheaf splits (i.e. direct sum of line bundles) in terms of the geometry of the singular scheme of the distribution.

The projective space case has been already considered. Indeed,
when $\iota_{X} = 4,$ L. Giraldo and A. J. Pan-Collantes showed in \cite{GP},
that the tangent sheaf of a foliation of dimension $2$ on $\mathbb{P}^3,$ splits if and only if its singular scheme $Z$ is aCM \cite[Theorem 3.3]{GP}.

\medskip

Let us now consider the quadric hypersurface $Q^3$ in $\p^4,$ i.e. $\iota_{X} = 3.$
Let $E$ a bundle on $Q^3.$ It is well known that if $E$ splits on $Q^3,$ then:
$$H^i(Q^3, E(t)) = 0  \;\; \mbox{for} \;\; 0<i<3, \;\; \mbox{for all} \;\; t \in \Z.$$

The following result characterizes when the tangent sheaf of a distribution of
dimension $2$ on $Q^3,$ is split or spinor.

\begin{thm}
Let $\sF$ be a distribution on $Q^3$ of codimension one,
such that the tangent sheaf $T_{\sF}$ is locally free.
If $T_{\sF}$ either splits as a sum of line bundles or is a spinor bundle, then $Z$ is arithmetically Buchsbaum, with $h^1(Q^3, I_{Z} (r-2)) = 1$ being the only nonzero intermediate cohomology for $H^i(I_{Z}).$
Conversely, if $Z$ is arithmetically Buchsbaum
with $h^1(Q^3, I_{Z} (r-2)) = 1$ being the only nonzero intermediate cohomology for $H^i(I_{Z})$ and $h^2(T_{\sF}(-2)) = h^2 (T_{\sF} (-1-c_1(T_{\sF}))) = 0,$ then $T_{\sF}$ either split or is a spinor bundle.
\end{thm}

\begin{dem}
Suppose that $T_\sF$ either split or
is a spinor bundle.
By considering the short exact sequence (\ref{eq:Dist Fano}), after tensoring with $\mathcal{O}_{Q^3}(t),$
and taking the long exact sequence of cohomology we get:
\begin{eqnarray} \label{cohom}
\cdots \rightarrow  H^1(Q^3,T_{\sF}(t)) \rightarrow H^1(Q^3, T Q^3(t)) \rightarrow H^1(Q^3, \mathcal{I}_{Z} (r+t)) \rightarrow \\ \nonumber
\rightarrow H^2(Q^3,T_{\sF}(t)) \rightarrow H^2(Q^3, T Q^3(t)) \rightarrow H^2(Q^3, \mathcal{I}_{Z} (r+t))\rightarrow \cdots
\end{eqnarray}

Since $T_{\sF}$ has no intermediate cohomology, we have that
$$H^1(Q^3,T_{\sF}(t)) = H^2(Q^3,T_{\sF}(t)) = 0  \; \mbox{for all} \; t \in \Z,$$
thus,
$$H^1(Q^3, T Q^3(t)) \simeq H^1(Q^3, \mathcal{I}_{Z} (r+t)).$$

By Bott's formula for quadrics (Theorem \ref{BottQ}), we have that $H^1(Q^3, T Q^3(t)) = 0$ for all $t \neq -2,$ i.e. $H^1(Q^3, T Q^3 (-2)) \neq 0.$ Hence, $H^1(Q^3, \mathcal{I}_{Z} (r-2)) \neq 0$ and $Z$ is arithmetically Buchsbaum.

To prove the converse, suppose that $Z$ satisfies the properties in the statement
and that $h^2(T_{\sF}(-2)) = h^2 (T_{\sF} (-1-c_1 (T_{\sF}))) = 0.$


Consider the long exact cohomology sequence (\ref{cohom}) for all $t \neq -2.$
By Theorem \ref{BottQ}, we obtain that $ H^1 (Q^3, TQ^3 (t)) = 0,$ for all $t \neq -2.$ Applying Serre duality and Bott's formula, respectively, we get $H^0(Q^3, \mathcal{I}_{Z} (r+t)) = H^3(\mathcal{O}_{Q^3}(-t-r-3)) = 0,$ for all $r \neq 2.$ Thus, we have $H^1(Q^3,T_{\sF}(t))=0.$

By Serre duality, we conclude that
$$0 = H^1(Q^3,T_{\sF}(t)) = H^2 (Q^3,T_{\sF}(s)), \; \mbox{for all}\; s \neq -1-c_1 (T_{\sF}),$$
where $s = -t-3-c_1 (T_{\sF}).$
As by hypothesis $h^2 (T_{\sF} (-1-c_1 (T_{\sF}))) = 0,$ we conclude that $H^2 (Q^3,T_{\sF}(s))=0$ for all $s.$ Then, for all $t \neq -2, \;T_{\sF}$ either split or is a spinor bundle.

Now, consider the long exact cohomology sequence (\ref{cohom}), with $t=-2.$
By Bott's formula, we have that $ H^0 (Q^3, TQ^3 (-2)) = 0$ and $ H^2 (Q^3, TQ^3 (-2)) = 0.$ By hypothesis, we get
$h^2(T_{\sF}(-2))=0$ and $h^1(Q^3, I_{Z} (r-2)) = 1.$
Moreover, applying Serre duality and Bott's formula, respectively, we get $H^0(Q^3, \mathcal{I}_{Z} (r-2)) = H^3 (\mathcal{O}_{Q^3}(-1-r)) = 0,$ for all $r \neq 2.$ So, from the exact sequence
$$ 0 \rightarrow H^1(Q^3,T_{\sF}(-2)) \rightarrow H^1(Q^3, T Q^3(-2)) \simeq \C \xrightarrow{\beta} H^1(Q^3, \mathcal{I}_{Z} (r-2)) \simeq \C  \rightarrow 0$$
we conclude that $H^1(Q^3,T_{\sF}(-2)) = 0,$ since $\beta$ is injective and $ker (\beta) = H^1(Q^3,T_{\sF}(-2)).$ Then, for $t=-2, \; T_{\sF}$ either split or is a spinor bundle.
Therefore, $T_{\sF}$ either split or is a spinor bundle, for all $t \in \Z.$

\end{dem}


Based on the result in \cite[Theorem 11.8]{MOJ}, we can construct examples of codimension one distribution on $Q^3.$


\begin{Example}
The bundle Spinor on $Q^3$ twisted by $\mathcal{O}_{Q^3}(1-t),$ is the tangent sheaf of
a codimension one distribution $\sF,$ for all $t \geq 0.$ Indeed, since the spinor bundle $S$ is globally generated, $S(t)$ is globally generated, for all $t \geq 0.$ 
Moreover, its rank $2$ reflexive sheaf on $Q^3,$ and
$S(t) \otimes TQ^3$ is also globally generated, since $TQ^3$ is globally generated. 
Now, apply \cite[Theorem 11.8]{MOJ} with $L =  \mathcal{O}_{Q^3}$ to obtain the desired codimension
one distribution
$$ 0 \rightarrow S^{\ast} (-t) \simeq S(1-t) \rightarrow TQ^3 \rightarrow \mathcal{I}_Z (r) \rightarrow 0.$$
\end{Example}

If $\iota_{X}=2,$ we characterize when the tangent sheaf of a distribution of
dimension $2$ on a smooth weighted projective complete intersection del Pezzo Fano threefold $X,$ has no intermediate cohomology. More precisely, we prove the following results:

\begin{Lemma} \label{coho-delPezzo}
Let $X$ be a smooth weighted projective complete intersection del Pezzo Fano threefold. Then, $H^2(X, \Omega^1_{X}(t)) = 0$ for $t > 4,$ and $H^1(X, \Omega^2_{X}(t)) = 0$ for $t > 10.$
\end{Lemma}

We omit the proof of the previous lemma, being a computation
of the cohomolology groups $H^p(X, \Omega_X^q (t))$ with $p,q \in \{1,2\}$ and $p \neq q,$ where $X$ is each one of the varieties with index two described in Section \ref{section:Fano3folds}. For this calculation, we use the sequences and cohomology formulas of the Section \ref{cohomologia}.
We summarize the computations in the following table:

\begin{table}[h] \label{tableH1} 
\center
\begin{tabular}{|c | c | c | c | c | c |  } \hline
$X$           &$ \iota_X $  & $H^2(X, \Omega^1_X(t))$   & $ t $     & $H^1(X, \Omega^2_X(t))$                    & $ t $ \\ \hline
$X_3$         &   2          & $0$                      &$t>1$           & $ 0 $                                      &  $t>4$            \\ \hline
$X_{2,2}$     &   2          &$0$                       &$t>0$           & $ 0 $                                      &  $t>2$  \\ \hline
$X_4$         &   2          &$0$                       &$t>2$           & $ 0$                                       &  $t>6$    \\ \hline
$X_6$         &   2          &$0$                       &$t>4$           & $ 0 $                                      &  $t>10$     \\ \hline
\end{tabular}
\caption{Values of $t$ for which $H^1(X,\Omega^2_{X}(t)) = H^2(X,\Omega^1_{X}(t)) = 0.$}
\end{table}

\begin{thm}
Let $\sF$ be a distribution of codimension one on a smooth weighted projective complete intersection del Pezzo Fano threefold $X,$ such that the tangent sheaf $T_{\sF}$ is locally free.
If $T_{\sF}$ has no intermediate cohomology, then $H^1(X, I_{Z} (r+t)) = 0$ for $t < -6$ and $t > 8.$
Conversely, if $H^1(X, I_{Z} (r+t)) = 0$ for $t < -6$ and $t > 8,$ and $H^2(X, T_{\sF}(t)) = 0$ for $t \leq 8$ and
$H^1(X, T_{\sF}(s))  = 0$ for $s \neq -t- \iota_X -c_1(T_{\sF}),$ then $T_{\sF}$ has no intermediate cohomology.
\end{thm}

\begin{dem}
Suppose that $T_\sF$ has no intermediate cohomology. By considering the short exact sequence (\ref{eq:Dist Fano}), after tensoring with $\mathcal{O}_{X}(t),$ and taking the long exact sequence of cohomology we get:
\begin{eqnarray} \label{cohomTX}
\cdots \rightarrow  H^1(X,T_{\sF}(t)) \rightarrow H^1(X, T X(t)) \rightarrow H^1(X, \mathcal{I}_{Z} (r+t)) \rightarrow \\ \nonumber
\rightarrow H^2(X,T_{\sF}(t)) \rightarrow H^2(X, T X(t)) \rightarrow H^2(X, \mathcal{I}_{Z} (r+t))\rightarrow \cdots
\end{eqnarray}
Since $T_{\sF}$ has no intermediate cohomology, we have that
$$H^1(X,T_{\sF}(t)) = H^2(X,T_{\sF}(t)) = 0  \; \mbox{for all} \; t \in \Z.$$
Thus, we get $H^1(X, T X(t)) \simeq H^1(X, \mathcal{I}_{Z} (r+t)).$

By Remark \ref{formula}, we have that $H^1(X, T X(t)) \simeq H^1(X, \Omega^2_{X} (t+2))$ and by Lemma \ref{coho-delPezzo}, we get $H^1(X, T X(t))= 0,$ for $t > 8.$ Moreover, by using Serre duality, we obtain $H^1(X, T X(t)) \simeq H^2(X, \Omega^1_{X} (-t-2))$ and thus, by Lemma \ref{coho-delPezzo}, we get $H^1(X, T X(t)) = 0,$ for $t < -6.$
Therefore, $H^1(X, \mathcal{I}_{Z} (r+t)) = 0$ for $t < -6$ and $t > 8.$

Conversely, suppose that $h^2(T_{\sF}(t)) = 0 $ for $t \leq 8.$ Consider the long exact cohomology sequence (\ref{cohomTX}). Applying Serre duality and Theorem \ref{flenner},
respectively, we get $H^2(X, TX(t)) = 0$ for $t > -2.$ Since by hypothesis,
$h^1(X, I_{Z} (r+t)) = 0$ for  $t < -6$ and $t > 8,$ we get $H^2(X, T_{\sF}(t)) = 0$ for $t > 8.$
Thus, $H^2(X, T_{\sF}(t)) = 0 $ for all $ t \in \Z.$

By Serre Duality and by Remark \ref{formula}, we conclude that
$$ 0 = H^2(X, T_{\sF}(t)) = H^1(X, T_{\sF}(s)),$$
where $s =  -t- \iota_X -c_1 (T_{\sF}).$  As by hipothesis $h^1 (T_{\sF} (s)) =0$ for $s \neq -t- \iota_X -c_1 (T_{\sF}),$ we conclude that $H^1(X, T_{\sF}(t)) = 0 $ for all $ t \in \Z.$
Therefore, $T_{\sF}$ has no intermediate cohomology.

\end{dem}

In \cite{AC}, Arrondo e Costa classify all aCM bundles on the Fano threefolds $X_d \subset P^{d+1}$ $d=3,4,5$
of index
$2.$ Such bundles are isomorphic to a twist of either $S_L,$ or $S_C,$ or
$S_E,$ where $L \subset X_d$ is a line, $C \subset X_d$ is a conic and $E \subset X_d$ is an elliptic curve of degree $d+2,$ respectively. Moreover, the authors showed that $S_L(1), S_C$ and $S_E(1)$ are globally generated vector bundles.

By using the result \cite[Theorem 11.8]{MOJ}, we can construct examples of codimension one distribution on $X_d.$
Consider $\mathcal{E} \simeq S_L(1); S_C; S_E(1).$ We have that $\mathcal{E}$ is generated by global sections.

\begin{Example}

Let $\mathcal{E}$ be a globally generated rank $2$ reflexive sheaf on $X_d.$ Then $\mathcal{E}^{\ast}(-t)$
is the tangent sheaf of a codimension one distribution  ${\sF},$ for all $t \geq 0.$ Indeed,
since $\mathcal{E}$ is globally generated, $\mathcal{E}(t)$ is globally generated, for all $t \geq 0$. Moreover, its rank $2$ reflexive sheaf on $X_d,$ and
$\mathcal{E}(t) \otimes TX_d$ is also globally generated, since $TX_d$ is globally generated.
Now, apply the Theorem \cite[Theorem 11.8]{MOJ} with $L =  \mathcal{O}_{X_d}$ to obtain the desired codimension
one distribution
$$ 0 \rightarrow \mathcal{E}^{\ast} (-t)  \rightarrow TX_d \rightarrow \mathcal{I}_Z (r) \rightarrow 0.$$

\end{Example}

\medskip

Finally, if $\iota_{X} = 1,$ we characterize when the tangent sheaf of a distribution of
dimension $2$ on a smooth weighted projective complete intersection prime Fano threefold $X,$ has no intermediate cohomology. More precisely, we prove the following results.

\begin{Lemma} \label{coho-prime}
Let $X$ be a smooth weighted projective complete intersection prime Fano threefold. Then, $H^1(X, \Omega^2_{X}(t)) = 0$ for $t > 5,$ and $H^2(X, \Omega^1_{X}(t)) = 0$ for $t > 3.$
\end{Lemma}

We will omit the proof of the previuos lemma, being a very similar computation
of the del Pezzo Fano's case.
We summarize the computations in the following table:

\begin{table}[h] \label{tableH2} 
\center
\hspace{-1.5cm}
\begin{tabular}{|c | c | c | p{4.2cm} | c | p{4.2cm} |  } \hline
$X$           &$ \iota_X $  & $H^2(X, \Omega^1_X(t))$   & $ \> \> \> \> t $     & $H^1(X, \Omega^2_X(t))$                    & $ \>\>\> \>t $ \\ \hline
$X_{2,3}$     &   1          &$0$                       &$t>2$ if $X_{2,3} \subset X_2$ and $t>3$ if $X_{2,3} \subset X_3$          & $ 0 $                                      & $t>5$ if $X_{2,3} \subset X_2$ and $t>5$ if $X_{2,3} \subset X_3$  \\ \hline
$X_{2,2,2}$   &   1          &$0$                        &$t>3$           & $ 0 $                                      & $t>5$  \\ \hline
$Y$           &   1           &$0$                      &$t>5$ if $Y \subset C_2$ and $t>5$ if $Y \subset X_4$         & $ 0 $                                      & $t>9$ if $Y \subset C_2$ and $t>3$ if $Y \subset X_4$  \\ \hline
\end{tabular}
\caption{Values of $t$ for which $H^1(X,\Omega^2_{X}(t)) = H^2(X,\Omega^1_{X}(t)) = 0.$}
\end{table}

\begin{thm}
Let $\sF$ be a distribution of codimension one on a smooth weighted projective complete intersection prime Fano threefold $X,$ such that the tangent sheaf $T_{\sF}$ is locally free.
If $T_{\sF}$ has no intermediate cohomology, then $H^1(X, I_{Z} (r+t)) = 0$ for $t < -4$ and $t > 4.$
Conversely, if $H^1(X, I_{Z} (r+t)) = 0$ for $t < -4$ and $t > 4,$ and $H^2(X, T_{\sF}(t)) = 0$ for $t \leq 4$ and
$ H^1(X, T_{\sF}(s))  = 0$ for $s \neq -t- \iota_X -c_1(T_{\sF}),$ then $T_{\sF}$ has no intermediate cohomology.
\end{thm}

\begin{dem}
Suppose that $T_\sF$ has no intermediate cohomology. Consider, for each $t \in \Z$, the exact sequence
$$0 \to T_{\sF}(t) \to TX(t) \to \mathcal{I}_Z(r+t) \to 0$$ and the long exact sequence of cohomology (\ref{cohomTX}).

Since $T_{\sF}$ has no intermediate cohomology, we get
$$H^1(X,T_{\sF}(t)) = H^2(X,T_{\sF}(t)) = 0 \; \mbox{for all} \; t \in \Z.$$
Thus, we have $H^1(X, T X(t)) \simeq H^1(X, \mathcal{I}_{Z} (r+t)).$

By Remark \ref{formula}, we get that $H^1(X, T X(t)) \simeq H^1(X, \Omega^2_{X} (t+1))$ and by Lemma \ref{coho-prime}, we get $H^1(X, T X(t)) = 0,$ for $t > 4.$ Moreover, using Serre duality, we obtain that $H^1(X, T X(t)) \simeq H^2(X, \Omega^1_{X} (-t-1))$ and thus, by Lemma \ref{coho-prime}, we get $H^1(X, T X(t)) = 0,$ for $t < -4.$
Therefore, $H^1(X, \mathcal{I}_{Z} (r+t)) = 0$ for $t < -4$ and $t > 4.$

Conversely, suppose that $h^2(T_{\sF}(t)) = 0 $ for $t \leq 4.$ Consider the long exact cohomology sequence (\ref{cohomTX}). Applying Serre duality and Theorem \ref{flenner},
respectively, we get $H^2(X, TX(t)) = 0$ for $t > -1.$ Since by hypothesis,
$h^1(X, I_{Z} (r+t)) = 0$ for  $t < -4$ and $t > 4,$ we get $H^2(X, T_{\sF}(t)) = 0$ for $t > 4.$
Thus, we have $H^2(X, T_{\sF}(t)) = 0 $ for all $ t \in \Z.$

By Serre Duality and by Remark \ref{formula}, we obtain that
$$ 0 = H^2(X, T_{\sF}(t)) = H^1(X, T_{\sF}(s)),$$
where $s =-t- \iota_X-c_1 (T_{\sF}).$  As by hipothesis $h^1 (T_{\sF}) (s) =0$ for $s \neq -t- \iota_X-c_1 (T_{\sF}),$ we conclude that $H^1(X, T_{\sF}(t)) = 0$ for all $t \in \Z.$
Therefore, $T_{\sF}$ has no intermediate cohomology.

\end{dem}

In \cite{Madonna}, Madonna classify aCM bundles on a prime Fano threefold $X_{2g-2} :=X$ of genus $g$. Let $l$ be a line in $X$ and consider that $\mathcal{E}$ is associated to a line $l$ in $X$ with $c_1=-1, c_2=1.$ It holds that $\mathcal{E}$ is generated by its global sections.
Again, by using the result \cite[Theorem 11.8]{MOJ}, we can construct examples of codimension one distribution on $X.$

\begin{Example}

Let $\mathcal{E}$ be a globally generated rank $2$ reflexive sheaf on a prime Fano threefold $X.$ Then $\mathcal{E}^{\ast}(-t)$
is the tangent sheaf of a codimension one distribution  ${\sF},$ for all $t \geq 0.$

\end{Example}


\section{Conormal sheaf vs. singular scheme} \label{cotangent sheaf}


In this section, all the distributions that we will consider will be defined for Fano threefolds, of codimension two, and with singular set of codimension $\geq 2.$

Looking at Definition \ref{distributionTF} we can alternatively define a foliation through a coherent subsheaf $N^{\ast}_{\sF}$ of $\Omega^1_X$ such that
\begin{enumerate}
  \item $N^{\ast}_{\sF}$ is integrable $(d N^{\ast}_{\sF} \subset N^{\ast}_{\sF} \wedge \Omega^1_X)$ and
  \item the quocient $\Omega^1_X / N^{\ast}_{\sF}$ is torsion free.
\end{enumerate}
The codimension of $\sF$ is the generic rank of $N^{\ast}_{\sF}$.

We characterize when the conormal sheaf is split, in terms of the geometry of the singular scheme of the distribution.
M. Corr\^ea, M. Jardim and R. Vidal Martins, showed
in \cite{CJV} that the conormal sheaf of a foliation of dimension one
on $\mathbb{P}^n$ splits if and only if its singular scheme is arithmetically Buchsbaum
with $h^1(\mathcal{I}_Z(d-1)) = 1$ being the only nonzero intermediate cohomology
\cite[Theorem 5.2]{CJV}.

Based in this result, we prove the following:

\begin{thm}
Let $\sF$ be an one-dimensional distribution on a smooth weighted projective complete intersection Fano threefold $X.$
 If $N^{\ast}_{\sF}$ is aCM, then $Z$ is arithmetically Buchsbaum, with $h^1(X, I_{Z} (r)) =1$ being the only nonzero intermediate cohomology for $H^i(I_{Z}).$
\end{thm}

\begin{dem}
Suppose that $N^{\ast}_{\sF}$ is aCM. For the case $\iota_X = 4,$ i.e. $X \simeq \mathbb{P}^3$, the result follows from Theorem $5.2$ in \cite{CJV}.

Consider, for each $t \in \Z$, the exact sequence
\stepcounter{thm}
\begin{equation} \label{exataX}
0 \rightarrow N^{\ast}_{\sF}(t) \rightarrow \Omega^1_{X}(t) \rightarrow \mathcal{I}_{Z} (r+t) \rightarrow 0,
\end{equation}
where $r$ is a integer such that $r = c_1(\Omega^1_{X}) - c_1(N^{\ast}_{\sF}).$
Taking the long exact sequence of cohomology we get:
\stepcounter{thm}
 \begin{eqnarray} \label{cohom-cot}
\cdots \rightarrow  H^1(X,N^{\ast}_{\sF}(t)) \rightarrow H^1(X, \Omega^1_{X}(t)) \rightarrow H^1(X, \mathcal{I}_{Z} (r+t)) \rightarrow \\ \nonumber
\rightarrow H^2(X,N^{\ast}_{\sF}(t)) \rightarrow H^2(X, \Omega^1_{X}(t)) \rightarrow H^2(X, \mathcal{I}_{Z} (r+t))\rightarrow \cdots
\end{eqnarray}
Since $N^{\ast}_{\sF}$ has no intermediate cohomology, we have that
$$H^1(X,N^{\ast}_{\sF}(t)) = H^2(X,N^{\ast}_{\sF}(t)) = 0, \; \mbox{for all} \; t \in \Z.$$
Thus, we get $H^1(X, \Omega^1_{X}(t)) \simeq H^1(X, \mathcal{I}_{Z} (r+t)).$
By Theorem \ref{flenner}, $H^1(X, \Omega^1_{X}(t))=0$ for all $t \neq 0,$ and $h^1(X, \Omega^1_{X})=1.$ Then,
we have $H^1(X,\mathcal{I}_{Z} (r)) \neq 0$ and $Z$ is arithmetically Buchsbaum.
\end{dem}

\medskip

If $X$ is a Fano threefold, meaning that $K^{-1}_X = \bigwedge^3 TX$ is ample, then the \textit{Kodaira vanishing Theorem} shows that $H^q(X, \mathcal{O}_X) = 0$ and $H^q(X, \mathcal{O}(K^{-1}_X)) = 0$  for $q>0,$ see \cite{Barron}.

A similar result for converse of Theorem $5.2$ in \cite{CJV} is the following:

\begin{thm}
Let $\sF$ be an one-dimensional distribution on a smooth weighted projective complete intersection Fano threefold $X,$ with index $\iota_X \in \{1,2,3\}.$
If $Z$ is arithmetically Buchsbaum
with $h^1(X, I_{Z} (r)) = 1$ being the only nonzero intermediate cohomology for $H^i(I_{Z}),$ and $h^2(N^{\ast}_{\sF}) = h^2 (N^{\ast}_{\sF} (-c_1 (N^{\ast}_{\sF}) - \iota_X)) = 0,$ then $N^{\ast}_{\sF}$ is aCM.
\end{thm}

\begin{dem}
Suppose that $h^2(N^{\ast}_{\sF}) = h^2 (N^{\ast}_{\sF} (-c_1 (N^{\ast}_{\sF}) - \iota_X)) = 0$ and that $Z$ is arithmetically Buchsbaum with $h^1(Q^3, I_{Z} (r)) = 1$ being the only nonzero intermediate cohomology.

\noindent Consider the long exact cohomology sequence (\ref{cohom-cot}), for all $t \neq 0.$
By Theorem \ref{flenner}, we have that $ H^1 (X, \Omega_{X}^1 (t)) = 0,$ for all $t \neq 0.$
Applying Serre duality and the Kodaira vanishing Theorem, respectively, we get
$$H^0(X, \mathcal{I}_{Z} (r+t)) = H^3(\mathcal{O}_{X}(-r-t- \iota_X))=0 \; \mbox{for} \; r \neq -\iota_X.$$
Thus, we get $H^1(X,N^{\ast}_{\sF}(t))=0$ for $t \neq 0.$

By Serre duality and by Remark \ref{formula}, we conclude that
$$0 = H^1(X,N^{\ast}_{\sF}(t)) = H^2 (X,(N^{\ast}_{\sF})^{\ast}(-t-\iota_X))= H^2 (X, N^{\ast}_{\sF}(-t-\iota_X-c_1 (N^{\ast}_{\sF}))).$$
Let $s = -c_1 (N^{\ast}_{\sF}) -t-\iota_X$ and $t \neq 0.$ Thus, we get
$$H^2 (X, N^{\ast}_{\sF}(s)) = 0 \; \mbox{for}\; s \neq - \iota_X -c_1 (N^{\ast}_{\sF}).$$
Since by hypothesis $h^2 (N^{\ast}_{\sF}(s)) = 0$ for $s =-\iota_X -c_1 (N^{\ast}_{\sF}),$  we conclude that \\ $H^2 (X,N^{\ast}_{\sF}(t))=0.$ Then, for all  $t \neq 0, \;N^{\ast}_{\sF}$ is aCM.

Now, consider the following piece of the long exact cohomology sequence (\ref{cohom-cot}), for $ t = 0:$
$$\ldots \rightarrow H^1(X,N^{\ast}_{\sF}) \rightarrow   H^1 (X, \Omega_{X}^1) \simeq \C \xrightarrow{\beta}
H^1(X, \mathcal{I}_{Z} (r)) \simeq \C \rightarrow 0.$$
The map $\beta$ is surjective and injective. Thus, $ker(\beta) = H^1(X,N^{\ast}_{\sF})$ is trivial, i.e.
\\ $H^1(X,N^{\ast}_{\sF}) = 0.$ As by hypothesis, $h^2(N^{\ast}_{\sF}) = 0,$ we conclude that
for $t=0, \; N^{\ast}_{\sF}$ is aCM.
Therefore, $N^{\ast}_{\sF}$ is aCM for all $t \in \Z.$

\end{dem}


\section{Properties of the singular locus of distributions} \label{connect}


In this section we analyze the properties of their singular schemes of codimension one holomorphic distributions on a threefold $X.$

\stepcounter{thm}
\subsection{Numerical Invariants}

Let $\sF$ be a codimension one distribution on threefold $X$ given as in the exact sequence (\ref{eq:Dist Fano}), with tangent sheaf $T_\sF$ and singular scheme $Z$.

Following \cite[section 2.1]{MOJ}, let $\mathcal{Q}$ be the maximal subsheaf of $\mathcal{O}_{Z/X}$ of codimension $>2$, so that one has an exact sequence of the form
\stepcounter{thm}
\begin{equation}\label{OZ->OC}
0 \to \mathcal{Q} \to \mathcal{O}_{Z/X} \to \mathcal{O}_{C/X} \to 0
\end{equation}
where $C\subset X$ is a (possibly empty) subscheme of pure codimension $2.$

The quotient sheaf is the structure sheaf of a subscheme $C\subset Z\subset X$ of pure dimension 1.

\begin{Definition}
If $Z$ is a 1-dimensional subscheme, then $Z$ has a maximal pure dimension 1
subscheme $C$ defining a sequence
\stepcounter{thm}
\begin{equation} \label{sequence I}
 0 \rightarrow \mathcal{I}_{Z} \rightarrow \mathcal{I}_{C} \rightarrow \mathcal{Q} \rightarrow 0,
\end{equation}
where $\mathcal{Q}$ is the maximal 0-dimensional subsheaf of $\mathcal{O}_{Z}.$
\end{Definition}

If $X$ is a Kähler manifold of dimension $n,$ and $Z \subset X$ is an analytic subset of codimension $k,$ then
\stepcounter{thm}
\begin{equation} \label{chern-I}
c_{k}(\mathcal{I}_{Z}) = (-1)^k (k-1)![Z],
\end{equation}

(for more details, see \cite{fulton}).

In the next theorem we obtain the following expressions for the Chern class of the tangent sheaf in terms of the cycle of $C.$ 
We will denote  $H=c_1(\mathcal{O}_X(1))$.

\begin{thm} \label{invariantes}
Let $\sF$ be a codimension one distribution on a threefold $X,$
with $\rho (X) = 1,$ given as in the exact sequence (\ref{eq:Dist Fano}), with tangent sheaf $T_\sF$ and singular scheme $Z$. Then,
$$c_2(T_\sF) = c_2(TX) - r \cdot K^{-1}_X + r^2H^2  -[C],$$ and
$$c_3(T_\sF) = c_3 (TX) - c_3(I_{Z/X})+ r H\cdot [C] - K^{-1}_X \cdot [C] -r H \cdot c_2(TX) + r^2 H^2   \cdot K^{-1}_X - H^3r^3.$$

\end{thm}

\begin{dem}
Considering the exact sequence (\ref{eq:Dist Fano}), we use that $c(TX)=c(T_\sF)\cdot c(I_{Z/X}(r))$ to obtain
\stepcounter{thm}
\begin{equation}\label{eq:IT}
\begin{array}{lll}
c_1(TX) & = & c_1(T_\sF)+c_1(I_{Z/X}(r)),\\
c_2(TX) & = & c_1(T_\sF)\cdot c_1(I_{Z/X}(r))+c_2(T_\sF)+c_2(I_{Z/X}(r)),\\
c_3 (TX) &= & c_3(T_\sF) + c_3(I_{Z/X}(r)) + c_1(T_\sF)\cdot c_2(I_{Z/X}(r)) + c_2(T_\sF)\cdot c_1(I_{Z/X}(r)).
\end{array}
\end{equation}

The first equation gives $c_1(T_\sF)=c_1(TX) - rH.$ From the exact sequence (\ref{sequence I}), it follows that $c_2(I_{Z/X}(r))=c_2(I_{C/X}(r))=[C]$, thus substitution into the second equation yields
$$ c_2(T_\sF)= c_2(T_X) - r H \cdot K^{-1}_X + H^2r^2 -[C]. $$

Moreover, the substituting the expressions for the first and second Chern classes into the third equation we obtain
\stepcounter{thm}
\begin{equation} \label{c3}
c_3 (TX) =  c_3(T_\sF) + c_3(I_{Z/X}(r)) + K^{-1}_X \cdot [C] -2r H\cdot [C] +r \cdot c_2(TX)- r^2H^2 \cdot K^{-1}_X + H^3r^3.
\end{equation}

Note that
\stepcounter{thm}
\begin{equation}\label{ast}
c_3(I_{Z/X}(r))=c_3(I_{Z/X})+r H\cdot c_2 (I_{Z/X})+ H^2r^3,
\end{equation}
while
\stepcounter{thm}
\begin{equation}\label{ast1}
c_2(I_{Z/X})= [C] - H^2r^2.
\end{equation}
Substituting \ref{ast1} into the equation \ref{ast}, we obtain
\stepcounter{thm}
\begin{equation}\label{ast2}
c_3(I_{Z/X}(r))=c_3(I_{Z/X})+ rH \cdot [C],
\end{equation}
and thus
\stepcounter{thm}
\begin{equation}\label{ast3}
c_3(T_\sF) = c_3 (TX) - c_3(I_{Z/X})+ rH \cdot [C] - K^{-1}_X \cdot [C] -rH \cdot c_2(TX) + r^2 H^2 \cdot K^{-1}_X - H^3r^3.
\end{equation}

\end{dem}

\stepcounter{thm}
\subsection{Connectedness}

Now, we will see when the singular locus of a codimension one
distribution is connected.
In \cite{MOJ}, O. Calvo-Andrade, M. Corr\^ea and M. Jardim present a homological
criterion for connectedness of the singular scheme of codimension $1$ distributions
on $\p^3.$

We will extend this criterion for others smooth weighted projective complete intersection Fano threefold with Picard number one.
Before, we need prove the vanishing of $H^1 (TX(-r))$ and $H^2 (TX(-r)).$

\begin{Lemma}
Let $X$ be a smooth weighted projective complete intersection Fano threefold with Picard number $\rho(X) = 1.$ Then, $H^1 (TX(-r)) = 0$ for $r > 6$ and $H^2 (TX(-r)) = 0$ for $r \neq \iota_X.$
\end{Lemma}

\begin{dem}
Suppose $X$ a smooth weighted projective complete intersection
Fano threefold.
If $\iota_X = 4,$ i.e. $X \simeq \p^3,$ by classical
Bott's formula we have that $H^1 (TX(-r)) = 0$ for all $r$ and $H^2 (TX(-r)) = 0$ for
$r \neq 2.$
If $\iota_X = 3,$ i.e. $X \simeq Q^3,$ by
Bott's formula for quadric, we have that $H^1 (TX(-r)) = 0$ for $r \neq 2$ and
$H^2 (TX(-r)) = 0$ for $r \neq 3.$
If $\iota_X = 2,$ by using Serre duality we get that $H^1 (TX(-r)) = H^2(\Omega^1_X(r-2)).$
By Table \ref{tableH1}, comparing the values of $t$ for which $H^2(\Omega^1_X(t)) = 0$ with $\iota_X = 2,$ we can see that the common vanishing of cohomology group for these varieties, occurs when $t > 4.$ Then, $H^2(\Omega^1_X(r-2)) = 0$ for $r > 6$ and $H^2 (TX(-r)) = 0$ for $r \neq 2$ by theorem \ref{flenner}.
If $\iota_X = 1,$ by using Serre duality we get $H^1 (TX(-r)) = H^2(\Omega^1_X(r-1)).$
By Table \ref{tableH2}, comparing the values of $t$ for which $H^2(\Omega^1_X(t)) = 0$ with $\iota_X = 1,$ we can see that the common vanishing of cohomology group for these varieties, occurs when $t > 3.$ Then, $H^2(\Omega^1_X(r-1)) = 0$ for $r > 4$ and $H^2 (TX(-r)) = 0$ for $r \neq 1$ by theorem \ref{flenner}.
Now, comparing the values of $r$ for which $H^1 (TX(-r)) = 0,$ we can see that the common vanishing of cohomology group considering all indices of $X,$  occurs when $r > 6$ and $H^2 (TX(-r)) = 0$ for $r \neq \iota_X.$

\end{dem}

\begin{thm}\label{con}
Let $\sF$ be a codimension one distribution with singular scheme $Z$ and let $X$ be a smooth weighted projective complete intersection Fano threefold with Picard number $\rho(X) = 1.$ If $h^2(T_\sF(-r))=0$ and $C \subset X, \; C  \neq \emptyset,$ then $Z$ is connected and of pure dimension 1, so that $T_\sF$ is locally free.
Conversely, for $r \neq \iota_{X}$, if $Z = C$ is connected, then $T_\sF$ is locally free and $h^2(T_\sF(-r))=0$.
\end{thm}

\begin{dem}
Twisting the exact sequence (\ref{eq:Dist Fano}) by $\mathcal{O}_{X}(-r)$ and passing to cohomology we obtain,
$$ H^1 (TX(-r)) \to H^1(I_{Z/X}) \to  H^2(T_\sF(-r)) \to H^2(TX(-r)). $$
By Lemma above, we get that $ H^1 (TX(-r))=0$ for $r > 6.$

If $h^2(T_\sF(-r))=0$, then $h^1(I_{Z/X})=0,$ for $r>6.$ It follows from the sequence
$$ 0 \to I_{Z/X} \to \mathcal{O}_{X} \to \mathcal{O}_{Z/X} \to 0 $$
that
$$ H^0(\mathcal{O}_{X} )\rightarrow H^0(\mathcal{O}_{Z/X})\rightarrow 0, $$
hence $h^0(\mathcal{O}_{Z/X})=1$. From the sequence (\ref{OZ->OC}), we get
$$ 0 \to H^0(\mathcal{Q}) \to H^0(\mathcal{O}_{Z/X}) \to H^0(\mathcal{O}_{C/X}) \to 0 $$
Thus either $h^0(\mathcal{O}_{C/X})=1$, and $\mathcal{Q}=0$ and $C$ is connected, or
$\mathrm{length}(\mathcal{Q})=1$ and $C$ is empty. This second possibility is not valid because by
hipothesis $C \neq \emptyset.$
It follows that $Z=C$ must be connect and of pure dimension 1, and thus, by Lemma \ref{local-free-tg-sheaf}, $T_\sF$ is locally free.

Conversely, assume that $Z=C$ is connected. Thus $Z$ must be of pure dimension 1, and Lemma \ref{local-free-tg-sheaf} implies that $T_\sF$ is locally free. It also follows that $h^1(I_{Z/X})=0,$ using Serre duality and the formula \ref{owci}. Since $h^2(TX(-r))=0$ for $r \neq \iota_{X}$, we conclude that $h^2(T_\sF(-r))=0$, as desired.
\end{dem}


\begin{Corollary}
If $\sF$ is a codimension one distribution on $X$ whose tangent sheaf splits as a sum of line bundles, then its singular scheme is connected.
\end{Corollary}

\begin{dem}
Assuming that $T_\sF = \mathcal{O}_X (r_1)\oplus \mathcal{O}_X (r_2)$, then clearly $h^2(T_\sF(-r))=0$, where $r=r_1+r_2.$
The result follows from Theorem \ref{con}.
\end{dem}

\begin{Corollary}
Let $\sF$ be a codimension one distribution on $X$ with locally free tangent sheaf.
If $T_\sF^{\ast}$ is ample, then its singular scheme is connected.
\end{Corollary}

\begin{dem}
We have, by Serre duality,
$$ H^2(T_\sF(-r)) \simeq H^1(T_\sF^{\ast}(r)\otimes K_{X})=H^1(T_\sF^{\ast}(r-c_1(TX))\otimes \mathcal{O}_X(c_1(TX))\otimes K_{X}). $$
Observe that $T_\sF^{\ast}(r-c_1(TX))\otimes \mathcal{O}_X(c_1(TX))\otimes K_{X} = T_\sF^{\ast}\otimes \det(T_\sF^{\ast})\otimes \mathcal{O}_X(c_1(TX))\otimes K_{X}$; since $T_\sF^{\ast}$ and $\mathcal{O}_X(c_1(TX))$ are ample, then, by Griffiths Vanishing Theorem \cite{Griffiths},
we get
$$
h^2(T_\sF(-r)) = h^1(T_\sF^*\otimes \det(T_\sF^*)\otimes \mathcal{O}_X(c_1(TX))\otimes K_{X}) =0 .
$$
The result follows from Theorem \ref{con}.
\end{dem}


\section{Stability} \label{stability}


In this section we will prove that   If $T_{\sF}$ is locally free with $c_1(T_{\sF}) = 0,$
then $T_{\sF}$ is either split or stable. Firstly, we will prove this for quadrics.

\begin{thm} \label{Q^3}
Let $\sF$ be a codimension one distribution on $Q^3.$ If $c_1(T_{\sF}) = 0,$
then $T_{\sF}$ is either split or stable.
\end{thm}

\begin{dem}
Since $c_1(T_{\sF}) = 0,$  we have the sequence
$$\sF: 0 \to T_{\sF} \to TQ^3 \to \sI_{Z/Q^3} (3) \to 0.$$  Suppose that $\sF$ is non stable and non split. Let $\alpha \in H^0(T_{\sF})$ be a nonzero section. This section induces a section $\sigma$ on $TQ^3.$ Letting $\sI_{\Gamma}:=\coker\alpha,$ we have the following diagram after applying the Snake Lemma:

$$ \xymatrix{
&   & 0 \ar[d]  & 0 \ar[d]  &        & \\
&   & \mO_{Q^3} \ar@{=}[r]  \ar[d]^-{\alpha}  &   \mO_{Q^3}  \ar[d]^-{\sigma}   &        & \\
&0\ar[r]& T_\sF  \ar[r] \ar[d]   & TQ^3 \ar[r] \ar[d]   & \sI_Z(3)  \ar@{=}[d] \ar[r] &0 \\
&0\ar[r]& \sI_{\Gamma}  \ar[r] \ar[d]   & \G \ar[r] \ar[d]   &\sI_Z(3) \ar[r] &0 \\
& &   0   &  0    &   &
} $$

Dualizing the last line of diagram, we obtain
$$0 \to \mO_{Q^3}(-3) \to \sG^{\ast} \to \mO_{Q^3} \to \inext^1(\sI_Z(3), \mO_{Q^3}) \to \inext^1(\sG, \mO_{Q^3}) \to
\inext^1(\sI_{\Gamma}, \mO_{Q^3}) $$

$$\to \inext^2(\sI_Z(3), \mO_{Q^3}) \to 0.$$
By the sequence $0 \to \sI_Z \to \mO_{Q^3} \to \mO_Z \to 0$ \label{ideais} we have the long exact sequence
$$0 \to  \Hom(\mO_Z, \mO_{Q^3}) \to \Hom(\mO_{Q^3}, \mO_{Q^3}) \to \Hom(\sI_Z, \mO_{Q^3}) \to \inext^1(\mO_Z, \mO_{Q^3}) \to \cdots.$$
Because $Z$ is a curve of codimension $2$ we have $\inext^i(\mO_Z, \mO_{Q^3}) = 0$ for $i<2,$ and thus
$$ \mO_{Q^3} = \Hom(\mO_{Q^3}, \mO_{Q^3}) \simeq \Hom(\sI_Z, \mO_{Q^3}).$$

Dualizing the first line of diagram, we obtain
$$0 \to \mO_{Q^3}(-3) \to \Omega^1_{Q^3} \to \sF^{\ast} \to $$

$$ \to \inext^1(\sI_Z(3), \mO_{Q^3}) \to \inext^1(TQ^3, \mO_{Q^3}) \to \inext^1 (\sF,\mO_{Q^3}) \to$$

$$ \to \inext^2(\sI_{Z}(3), \mO_{Q^3}) \to \inext^2(TQ^3, \mO_{Q^3}) \to \cdots $$

A sheaf $E$ is locally free if and only if $\inext^i(E, \mO_X)=0$ for each $i \geq 1.$ Thus, $$\inext^1 (\sF,\mO_{Q^3}) = \inext^2(TQ^3, \mO_{Q^3}) = 0$$ and by the sequence above we get $\inext^2(\sI_{Z}(3), \mO_{Q^3})=0.$

\bigskip

Twisting the sequence \ref{ideais} by $\mO_{Q^3}(3)$ and dualizing, we obtain

$$ \cdots \to \inext^1(\mO_Z(3), \mO_{Q^3}) \to \inext^1(\mO_{Q^3}(3), \mO_{Q^3}) \to \inext^1(\sI_Z(3), \mO_{Q^3}) \to$$

$$\to \inext^2(\mO_Z(3), \mO_{Q^3}) \to \inext^2(\mO_{Q^3}(3), \mO_{Q^3}) \to \inext^1(\sI_Z(3), \mO_{Q^3}) \to$$

$$\to \inext^3(\mO_Z(3), \mO_{Q^3}) \to \inext^3(\mO_{Q^3}(3), \mO_{Q^3}) \to \cdots$$

We have $\inext^i(\mO_{Q^3}(3), \mO_{Q^3}) = 0 $ for each $i \geq 1$ and $\inext^2(\sI_{Z}(3), \mO_{Q^3})=0$ by previous computations. Therefore, $\inext^3(\mO_Z(3), \mO_{Q^3}) = 0$ and $$\inext^1(\sI_Z(3), \mO_{Q^3}) \simeq \inext^2(\mO_Z(3), \mO_{Q^3}) = \omega_Z(-3),$$ where $\omega_Z$ is the dualizing sheaf of the curve $Z.$

\bigskip

Similarly, from the sequence $0 \to  \sI_{\Gamma} \to \mO_{Q^3} \to \mO_{\Gamma} \to 0,$ we get $$\inext^1(\sI_{\Gamma}, \mO_{Q^3}) \simeq \inext^2(\mO_{\Gamma}, \mO_{Q^3}).$$ Since $\Gamma$ is a codimension two, we have the local fundamental isomorphism
$$\inext^2(\mO_{\Gamma}, \mO_{Q^3}) \simeq  \Hom( \det(\sI_{\Gamma}/\sI^2_{\Gamma}), \mO_{Q^3} ) = \mO^{\ast}_{Q^3} \otimes \mO_{Q^3} = \mO_{Q^3}.$$ Therefore, $\inext^1(\sI_{\Gamma}), \mO_{Q^3}) =  \mO_{Q^3}.$
Thus, we obtain $$\supp(\inext^1(\sF, \mO_{Q^3})) \subset \Sing (\sI_{\Gamma}) = (\alpha=0) \subset \supp(\inext^1(\sG, \mO_{Q^3})) = \Sing (\sG).$$

\bigskip

Since $\Gamma$ is smooth, the Chern classes of $T_\sF$ are given by $c_1(T_{\sF})=0$ and
 $c_2(T_{\sF})= \deg (\Gamma) = d.$

\bigskip

Since $\omega_{Q^3} \simeq \mO_{Q^3}(-3),$ we have  $\omega_{\Gamma} \simeq \mO_{\Gamma}(c_1-3).$
If $\Gamma$ is an irreducible nonsingular curve of genus $g,$ then $\omega_{\Gamma}$ is the canonical sheaf,
which has degree $2g-2.$ Therefore,
$$2g-2 = d(c_1-3) = c_2(c_1-3) = -3c_2 \label{genus}.$$

By using the Theorem \ref{invariantes}, we have $c_2(T_\sF) = 4H^2 - H\cdot [Z].$ Since $\Gamma$ is a curve, it has at least $\deg (\Gamma) = 1,$ thus $ H\cdot [Z]$ is at most three. Then, if $H\cdot [Z]=0 \;(\deg(\Gamma) = 4), $ we have $2g-2=-12 \Rightarrow g = -5$ and if $H\cdot [Z]=2\; (\deg(\Gamma) = 2),$ then $2g-2=-6 \Rightarrow g = -2.$

On the  other hand, since $T_{\sF}$ is locally free, it has $c_3(T_{\sF})=0.$ By Riemann-Roch theorem $c_3(I_{Z}) = (2g-2)H^3+c_1(Q^3)c_2(I_{Z}).$ Since $Z$ is a curve, we have $c_2(I_Z) = [Z] $ and we know that  $c_1(Q^3) = 3H.$ Then, $c_3(I_Z) = (2g-2)H^3+3H\cdot[Z].$ Replacing these values in the Theorem \ref{invariantes}, we have $$(2g-2)H^3=-10H^3-3H\cdot [Z]. \label{genus2}$$
Now, by using the formula \ref{genus2}, if $H\cdot [Z]=0,$ we have a curve of genus $g=-4$ and if $H\cdot [Z]=2,$ we have a curve of genus $g=-7.$ It is a contradiction. Therefore, $T_{\sF}$ is either stable or split.

\end{dem}

\begin{Remark} \label{isolated}
Note that if the subfoliation has only isolated singularities then $T_{\sF}$ has to be split. Indeed, $T_{\sF}$ would have a section that does not vanish on a curve, so it would have to be split.
\end{Remark}

\begin{thm}
Let $X$ be a Fano threefold with Picard number one and $\sF$ be a distribution of codimension 1. If $T_{\sF}$ is locally free with $c_1(T_{\sF}) = 0,$
then $T_{\sF}$ is either split or stable.
\end{thm}

\begin{dem}
If $\iota_X = 4,$ then the result  follows from  \cite[Theorem 9.5]{MOJ}.
If $\iota_X = 3,$ it follows from Theorem \ref{Q^3}.
If $\iota_X = 2,$ suppose that $T_{\sF}$ has a global section, so such section induce a global section of $TX.$ It follows from
[\cite{LPT}, Corollary 7.4] that the only Fano of index two with global field is $X_5.$ And by [\cite{LPT}, Lemma 7.2], every global field of $X_5$ has only isolated singularities. Thus, $T_{\sF}$ is split by Remark \ref{isolated}. It is a contradiction.
If $\iota_X = 1,$ follows from the proof of Theorem $8.1$ in \cite{LPT} that $X$ is birational to $X_5.$ As in $X_5,$ every global field of $X$ has only isolated singularities. Thus, $T_{\sF}$ is split by Remark \ref{isolated}. Once again this  is a contradiction.

\end{dem}

\end{document}